
\baselineskip=14pt
\parskip=10pt

\magnification=\magstephalf

\parindent=0pt
\overfullrule=0in
\def\frac#1#2{{#1 \over #2}}

\centerline
{ \bf
An Empirical Method for Solving (Rigorously!) Algebraic-Functional Equations of the Form
}
\centerline
{
 $ {\bf F( P(x,t),P(x,1),x,t) \equiv 0}$
}

\rm
\bigskip
\centerline
{\it By Ira M. GESSEL and Doron ZEILBERGER}

\bigskip

{\bf Preface: It is Never Too Late}

One of the  references in Doron Zeilberger's article [Z1], written more than $23$ years ago, 
was to an article, labeled {\it in planning}, with the authors and title of the present article.
We must have both forgotten our commitment. But {\it better late than never}, so with a delay of almost a quarter-century, at long
last we got around to writing the promised article. But it is {\it just as well\/} that we waited, since
now computers are much faster, and we are much better programmers!

{\bf Functional Equations}

Many enumeration problems reduce to solving  a functional equation for a generating function with respect
to one or more variables that we {\bf do} care about, and one or more variables, that we {\bf don't} care
about, called {\it catalytic variables}. At the {\it end of the day\/} we set all the catalytic variables to~$1$ (or sometimes $0$, and possibly other specific numbers).
Even though---at least initially---we do not care
about these extra variables (corresponding to auxiliary combinatorial statistics), they are needed
to set up the {\it combinatorial reaction\/}, so to speak, only to be discarded, like their
chemical namesakes, once the `reaction' is finished.

Initially humans solved these, one at a time, using {\it ad hoc\/} human ingenuity, and this goes back
to the pioneering work of Tutte and his school on counting maps, in the early 1960s.

More recently, the favorite method became the powerful and sophisticated {\it kernel method\/}
that had scored many triumphs in the hands of
such virtuosi as CNRS Silver-Medalist Mireille Bousquet-M\'elou and other people. 
For  lucid and engaging overviews, see the slides of the talks [B1] and [B2].
Alas, this method, in addition to requiring a lot of {\it human ingenuity}, is also
very {\bf human-labor-intensive}.

It turns out, that in {\it many cases\/} (perhaps all!), there is an alternative,
much simpler, approach, based on empirical {\it guessing}, {\bf yet it is fully rigorous!}
Of course, this method requires the help of our silicon colleagues.

{\bf The Zeilberger Gordian Knot}

In Zeilberger's proof [Z1] of Julian West's [W] conjectured explicit expression for
the number of $2$-stack-sortable permutations, it was necessary to solve the  functional equation
$$
f(x,t)= \frac{1}{1-xt} + \frac{xt (f(x,1)-t f(x,t))(f(x,1)-f(x,t))}{(1-t)^2}  .
\eqno(\rm FunEq)
$$
By clearing denominators, this equation can be written more compactly as
$$
F(\, f(x,t) \, , \, f(x,1)\, , \, x \, ,\, t) \, \equiv \, 0 ,
$$
for some polynomial $F$ of four variables.

The idea is extremely simple. Since $f(x,1)$ was conjectured to be a certain known {\it algebraic\/}  formal power series,
why not {\it guess\/} that this is also the case for the two-variate formal power series $f(x,t)$; i.e., let the
computer {\bf guess} a polynomial in three variables---let's call it $G$---such that
$$
G( f(x,t) \, , \, x \, , \, t) \, \equiv \, 0 .
\eqno(\rm AlgEq)
$$
Once guessed, it is purely routine to prove our guess rigorously. Since both $(\rm FunEq)$ and $(\rm AlgEq)$ have unique formal power series solutions,
after we define $g(x,t)$ to be the unique solution of $\rm (AlgEq)$, the verification that
$$
F(g(x,t) \, , \,  g(x,1) \, , \, x \, , \, t) \, \equiv \, 0  ,
$$
is a routine calculation in the {\it Sch\"utzenberger ansatz} [Z2]. In fact, since we know {\it a priori\/} that
the left side satisfies {\it some\/} algebraic equation, all we need is to  bound the degrees, and
check that the first few coefficients (in $x$) are identically $0$. Since that's how we got it in the first place, we already know that!
Nevertheless, the pedantic purist may want to bound the degrees exactly.

{\bf The Gessel Shortcut}

{\bf Alas}, guessing the three-variate polynomial $G$ takes a very long time.
After the first draft of [Z1] was written, Ira Gessel made the following observation, reproduced at
the very end of the final version of [Z1], that we now reproduce.

``{\bf  Epilogue: How The Proof Could Have Been Found}

July 2, 1991:
The first proof of any conjecture is seldom the shortest. It turns
out that the present proof is no exception. Ira Gessel made
the brilliant observation that steps 2--5 can be replaced by
the following.

{\bf Step $2'$}: Conjecture a polynomial $I$, of two variables, such that $I(P(x),x))=0$
where $P(x)=f(x,1)$. To prove it rigorously,
we must show that the unique $ f(x,t)$ that satisfies
$$F( f(x,t), f(x,1),x,t) \equiv 0$$ is such that $I( f(x,1),x) \equiv 0$.
Let's write,

$$
  {\rm(i)} \  F( f(x,t) , Q(x) , x, t ) \equiv 0,  \quad {\rm(ii)} \ f(x,1) = Q(x),\quad {\rm(iii)}\ I( Q(x) , x ) \equiv 0.
$$

We have to prove that (i)+(ii) implies (iii). But note that
(i)+(ii) have a unique solution, and (i)+(iii) have a unique
solution, and we must show that these are the same. So it's
enough to show that (i)+(iii) implies (ii). Taking the resultant
of $F$ and $I$ w.r.t. $Q(x)$ gives the algebraic equation
$G( f , x, t) \equiv 0$ found empirically, and very painfully, in step 3.
Proceeding as in step 5, we see that indeed $Q(x)=P(x)$.
This observation is the {\it leitmotif\/} of a paper that Ira Gessel
and I hope to write.''

\vfill\eject

{\bf Implementation}

This method turns out to be applicable to many other functional equations. Using Maple, it is
easy to follow Gessel's advice. It is very fast to guess a polynomial of two variables,
let's call it $I$, such that
$$
I(f(x,1),x) \, \equiv \, 0  .
$$
By the assumption, we have the relation
$$
F( f(x,t) \, , \, f(x,1) \, , \, x \, , \, t) \, \equiv \, 0  .
$$

Now eliminate $f(x,1)$ from both equations, using, say, Maple's built-in command {\tt resultant}, and get an algebraic
equation linking $f(x,t)$, $x$, and $t$, {\bf without} $f(x,1)$. Now replace $f(x,t)$ by $f(x,1)$, and $t$ by $1$, and get 
a polynomial in $f(x,1)$ and $x$ and make sure that it is a nonzero multiple of $I$. This is so much faster than
the original approach, and all the steps are fully automatic.

{\bf Linear Recurrence Equations with Polynomial Coefficients}

It is well-known, and easy to see (and implement, e.g., in the Maple package {\tt gfun} described in [SZ])
that once a formal power series satisfies an algebraic equation, as above, it also satisfies a
{\it linear differential equation with polynomial coefficients\/} (i.e., is $D$-finite), and hence
the enumerating sequence itself, satisfies a {\it linear recurrence equation with polynomial coefficients} (i.e., is $P$-recursive).
All this can be found automatically, and in fact, since we know {\it a priori\/} that such a recurrence is bound to exist,
it is completely legitimate to guess it empirically. If we are lucky, and that recurrence happens to be
{\it first-order}, then we get a closed-form `elegant' expression, like in the case of [W], first proved in [Z1].

{\bf The Maple Package FunEq}

Everything here is fully implemented in the Maple package {\tt FunEq}, available from the front
of this article: {\tt http://www.math.rutgers.edu/\~{}zeilberg/mamarim/mamarimhtml/funeq.html}.
It also contains quite a few sample input and output files that readers are welcome to extend.
In particular, we give fully rigorous automatic proofs of all the results in [CJS] (and many other ones, where
the answer is not `nice'), as well as a two-second proof of the main result of [Z1], and we solve numerous other functional equations
that we picked more or less at random, just to test the method.

{\bf Future Work: More Catalytic Variables; Beyond the Sch\"utzenberger Ansatz}

The present, naive, guessing approach should be
applicable to functional equations with more than one catalytic variable, but then,
according to [B1] and [B2], one may have to go to other ansatzes, first $D$-finite,
and then formal power series that satisfy an {\it algebraic differential equation\/} rather than a linear one.
Sooner or later, things would become too difficult even for computers, but one can do lots
of shortcut tricks, and we are sure that the present empirical-yet-rigorous approach should be
extendible, and {\it implementable}. Whether you would get nice results that humans can
enjoy remains to be seen, and is rather unlikely. Hence it may not be worth the effort.

{\bf References}

[B1] M. Bousquet-M\'elou, {\it Enumeration with ``catalytic'' parameters: a survey},
Three lectures at the 67th S\'eminaire Lotharingien de Combinatoire, Bertinoro, Italy, September 2011. Available from 
{\tt http://www.mat.univie.ac.at/\~{}slc/wpapers/s67vortrag/bousquet.pdf} \quad .

[B2] M. Bousquet-M\'elou, {\it M\'eli-m\'elo(u) de combinatoire}, 
expos\'e pour la remise de la m\'edaille d'argent du CNRS, Bordeaux, Oct. 2, 2014. Available from \hfill\break
{\tt http://www.labri.fr/Perso/\~{}bousquet/Exposes/medaille.pdf} \quad .

[CJS] R. Cori, B. Jacquard, and G. Schaeffer, {\it  Description trees for some families of planar maps},
``Proceedings of the 9th Conference on Formal Power Series and Algebraic Combinatorics '' (Vienna, 1997), 196--208. 
Available from  \hfill \break
{\tt http://www-igm.univ-mlv.fr/\~{}fpsac/FPSAC97/ARTICLES/Schaeffer.ps.gz} \quad .

[SZ] B. Salvy and P. Zimmermann, {\it GFUN: a Maple package
for the manipulation of generating and holonomic functions
in one variable}, ACM Trans. Math. Software {\bf 20} (1994), 163--177.

[W] J. West, {\it Sorting twice through a stack}, Theoretical Computer Science {\bf 117} (1993), 303--313.

[Z1] D. Zeilberger, 
{\it A proof of Julian West's conjecture that the number of two-stack-sortable permutations of length
$n$ is $2(3n)!/((n+1)!\,(2n+1)!)$}, Discrete Mathematics {\bf 102} (1992), 85--93. Available from
{\tt http://www.math.rutgers.edu/\~{}zeilberg/mamarimY/julian.pdf} \quad .

[Z2] D. Zeilberger,
{\it   An enquiry concerning human (and computer!) [mathematical] understanding},
in: ``Randomness and Complexity, From Leibniz to Chaitin'',   Cristian S. Calude, ed.,
World Scientific, Singapore, 2007, 383-410.  Available from \hfill\break
{\tt http://www.math.rutgers.edu/\~{}zeilberg/mamarim/mamarimhtml/enquiry.html} \quad .
\bigskip
\hrule
\bigskip

Ira M. Gessel,   Department of Mathematics, MS 050, Brandeis University, Waltham, MA 02453-2728, USA. 
\hfill \break
gessel at brandeis dot edu;  \quad {\tt http://people.brandeis.edu/\~{}gessel/} \quad .
\bigskip

\hrule

\bigskip
Doron Zeilberger, Department of Mathematics, Rutgers University (New Brunswick), Hill Center-Busch Campus, 110 Frelinghuysen
Rd., Piscataway, NJ 08854-8019, USA. \hfill \break
zeilberg at math dot rutgers dot edu;  \quad {\tt http://www.math.rutgers.edu/\~{}zeilberg/} \quad .

\bigskip
\hrule

\bigskip
EXCLUSIVELY PUBLISHED IN  The Personal Journal of Shalosh B. Ekhad and Doron Zeilberger ({ \tt http://www.math.rutgers.edu/\~{}zeilberg/pj.html}),
Ira Gessel's website,  and {\tt arxiv.org} \quad .
\bigskip
\hrule
\bigskip
{\bf Dec. 28, 2014}

\end